\documentstyle[12pt,leqno]{article}
\textheight9in         \def\DATE{January 7, 1998}
\newtheorem{theorem}{Theorem}[section]

\newtheorem{odstavec}[theorem]{}

\newtheorem{proposition}[theorem]{Proposition}

\catcode`\@=11

\def\@begintheorem#1#2{\it \trivlist \item[\hskip
 \labelsep{\bf #1\ #2.}]}
\def\@opargbegintheorem#1#2#3{\it \trivlist\item[\hskip%
 \labelsep{\bf #1\ #2.\ (#3)}]}
\def\@endtheorem{\endtrivlist}

\def\@listI{\leftmargin\leftmargini \parsep 1pt plus 2.5pt
 minus 1pt\topsep 5pt plus 2pt minus 3pt%
 \itemsep 0pt plus 2.5pt minus 1pt}
\let\@listi\@listI
\@listi

\def\@sect#1#2#3#4#5#6[#7]#8{\ifnum #2>\c@secnumdepth%
 \def \@svsec {}\else \refstepcounter {#1}\edef \@svsec%
 {\csname the#1\endcsname. \hskip .1em }\fi \@tempskipa%
 #5\relax \ifdim \@tempskipa >\z@ \begingroup #6\relax%
 \@hangfrom {\hskip #3\relax \@svsec }{\interlinepenalty%
 \@M #8.\par }\endgroup \csname #1mark\endcsname {#7}%
 \addcontentsline {toc}{#1}{\ifnum #2>\c@secnumdepth%
 \else \protect \numberline {\csname the#1\endcsname. }%
 \fi #7}\else \def \@svsechd {#6\hskip #3\@svsec #8.%
 \csname #1mark\endcsname {#7}\addcontentsline {toc}{#1}%
 {\ifnum #2>\c@secnumdepth \else \protect \numberline%
 {\csname the#1\endcsname. }\fi #7}}\fi \@xsect {#5}}

\def\section{\@startsection {section}{1}{\z@ }%
 {-3.5ex plus -1ex minus -.2ex}{2.3ex plus .2ex}{\bf }}

\def\thebibliography#1{%
 \section *{References.\@mkboth {REFERENCES}{REFERENCES}}%
 \list {[\arabic {enumi}]}{\settowidth \labelwidth {[#1]}%
 \leftmargin \labelwidth \advance \leftmargin \labelsep %
 \usecounter {enumi}} \def \newblock %
 {\hskip .11em plus .33em minus -.07em} \sloppy \clubpenalty 4000%
 \widowpenalty 4000 \sfcode`\.=1000\relax}

\def\@maketitle{%
 \newpage \null \vskip 2em
 \begin{center}
{\Large\bf \@title \par }
 \vskip 1.5em
 {\large \lineskip .5em
 \begin {tabular}[t]{c}\@author
 \end{tabular}\par} 
 \end{center}
  \vskip .8em}

\def\abstract{%
\if@twocolumn \section *{Abstract}
 \else \small\quotation\noindent{\bf Abstract.}\fi}

\def\ps@myheadings{\let\@mkboth\@gobbletwo
\def\@oddhead{\ifnum\count0=1 \hfill\else
\rightmark \hfil \rm \thepage\fi}%
\def\@oddfoot{\ifnum\count0=1 \hfill \rm 1 \hfill \else
\hfill\fi}
\def\@evenhead%
{\rm\leftmark\hfil\rm\thepage}%
\def\@evenfoot{}\def\sectionmark##1{}
\def\subsectionmark##1{}}

\catcode`\@=13

\def\G{{\sf G}} \def\F{{\cal F}}
\def\P{{\cal P}}

\def\bbb{\bf}        \def\Ass{{\it Ass}}
\def\SS{{\bbb S}}    \def\Comm{{\it Comm}}
                     \def\Lie{{\it Lie}}
    \def\T{{\sf T}}
\def\theor{{\bf T}}  \def\Q{{\cal Q}} 
\def\CC{{\bbb C}}    \def\bk{{\bf k}}
\def\id{{1\!\!1}}
\def\End{\mbox{\rm End}} \def\END{\mbox{\bf E}nd}
\def\Hom{\mbox{\rm Hom}}
\def\ot{\otimes}

\def\dual#1{{#1}^{\star}}
\def\vert#1{\mbox{\rm vert}(#1)}
\def\edg#1{\mbox{\rm edg}(#1)}

\def\leg#1{\mbox{\rm leg}(#1)}
\def\set#1{\{#1\}}
\def\cobar#1{{\Omega(\dual {#1})}}

\def\otexp#1#2{{#1}^{\otimes #2}}
\def\scalar#1#2{{\langle #1|#2\rangle}}

\def\ColoredTrees#1#2#3{{\cal T}^{#1}_{#2}(#3)}
\def\Graphs#1#2{{\cal G}_{#1}(#2)}
\def\ColoredGraphs#1#2#3{{\cal G}^{#1}_{#2}(#3)}
\def\ColoredGraphsPrime#1#2#3{{\cal G}'^{#1}_{#2}(#3)}

\begin{document}
\baselineskip18pt plus 2pt minus 1pt
\parskip3pt plus 1pt minus .5pt

\title{Cyclic operads and homology of graph complexes}

\author{Martin Markl%
\thanks{Supported by the grant GA AV \v CR \#1019507} 
\\
\normalsize 1st version: January 6, 1996; this version: January 7, 1998}

\maketitle
\begin{abstract}
We will consider $\P$-graph complexes, where $\P$ is a cyclic operad.
$\P$-graph complexes
are natural generalizations of Kontsevich's graph complexes -- for
$\P=\Ass$ it is the complex of ribbon graphs, for $\P= \Comm$ the
complex of all graphs. We construct a `universal 
class' in the cohomology of the graph complex with coefficients in a
theory. The Kontsevich-type invariant is then an evaluation, on a
concrete cyclic algebra, of this class. 

We also explain some results of
M.~Penkava and A.~Schwarz on the construction of an invariant from a
cyclic deformation of a cyclic algebra.
Our constructions are illustrated by a `toy model' of tree complexes.

\vskip 3mm
\noindent
\begin{tabular}{ll}%
\hskip-.55em{\bf Plan of the paper:}&%
\ref{intro}.~Introduction
\\
{}&\ref{Kolin}.~Warming up 
\\
{}&\ref{3ii}.~Graph complexes
\\
{}&\ref{4}.~The cycle
\\
{}&\ref{5}.~Lacunar graphs 
\\
{}&\ref{6}.~Appendix
\end{tabular}
\end{abstract}

\section{Introduction}
\label{intro}

In~\cite{kontsevich:91}, 
M.~Kontsevich constructed, for any cyclic $A_\infty$-algebra, an
element in the cohomology $H^*({\cal M}_{g,n};{\bf C})$ of the coarse moduli
space of smooth algebraic curves of genus $g$ with $n$ unlabeled punctures. 
His construction is based on a certain
combinatorial representation of ${\cal M}_{g,n}$ -- the graph complex -- 
and involves an
$A_\infty$-algebra. 
It resembles the state-sum model for the
Jones polynomial of a link
(in fact, it {\em is\/} a state sum). The aim of this note
is to give a conceptual understanding of the existence of these classes.

Because of the resemblance mentioned above, it would be helpful to
summarize the progress in the
understanding the quantum-group-type invariants of links.

\noindent
{\it \underline{1st} \underline{ste}p.}
The simplest state-sum model based on the canonical $R$-matrix of the 
quantum~$\mbox{\it sl\/}_2$. 

\noindent
{\it \underline{2nd} \underline{ste}p.}
A state-sum model related to the quantization of a general (semi-simple) Lie
algebra. 

\noindent
{\it \underline{3rd} \underline{ste}p.}
For a $\bk$-linear rigid braided monoidal category ${\cal C}$, 
each tangle $T$ (= `open' link) can be
interpreted as a morphism in ${\cal C}$, i.e.~as an element of ${\cal
C}(V^{\odot m},V^{\odot n})$, where $m$ (resp.~$n$) denotes the number of
input (output) strings of $T$.
A link is a closed tangle ($m = n = 0$) and we get an element of ${\cal
C}(1,1)$ ($1$ is the unit element of ${\cal C}$) which is a {\it
number\/}, because ${\cal C}(1,1) = \bk$. 

The construction of M.~Kontsevich mentioned above would correspond to the 
1st~step of the above imaginary 
list. To accomplish the
remaining two steps we need first a generalization of graph complexes to more
general types of algebras. This generalization was described, for Lie and
commutative associative algebras, by Kontsevich himself~\cite{kontsevich:93}, 
but `graph complexes' can be defined 
for algebras over an arbitrary {\it cyclic\/} (the
cyclicity is absolutely essential) operad. Such a generalization was, in
fact, given in~\cite{getzler-kapranov:preprint94} -- 
every cyclic operad can be naturally
considered as a modular operad, and the appropriate `graph complex' is the
Feynman transform of this operad introduced in the above mentioned paper. 

We will then show that there
exists a `universal cohomology class' (Proposition~\ref{central})
such that the invariant related to a
concrete algebra with an invariant scalar product
is a specialization (or evaluation) of this universal
class (\S\ref{cycle}, \S\ref{examples2}). 
The universality means that the class `contains' all
special invariants. The construction of this class was made possible by a
very explicit understanding of the structure of $\bk$-linear {\small PROP}s
or `theories' achieved in~\cite{markl:JPAA96,markl:ws93}. This class is not
only `universal' but also the `simplest possible' in the sense that it uses
only `generic' properties of objects.

As the first approximation of the understanding 
we offer the following comment.
We are going to
construct a complex and {\it simultaneously \/} a class in the homology
of this complex. The following analogy is helpful.
For any vector space $V$, the tensor product 
$\dual V \ot V$ of the dual $\dual V =
\Hom(V,\bk)$ 
and $V$ contains the `canonical element' $\eta \in \dual V \ot V$. If we
pick a basis $(e_i)_{i\in I}$ of $V$, 
we may give a `state-sum-type' definition of
$\eta$ as $\eta :=
\sum_{i\in I} \dual e_i \ot e_i$, where $(\dual e_i)_{i\in I}$ 
is the dual basis. A `categorical'
definition says just that $\eta$ corresponds, by duality, 
to the identity map $\id : V
\to V$. 

Observe that neither of the two `parts' (`$\dual V$'-part and
`$V$'-part) of $\eta$ can exists independently. 
In our analogy, the `$V$'-part is a
coloring of a graph, while the `$\dual V$'-part is the coefficient of
the cycle representing our canonical class, the coloring being given by an
element of our cyclic operad, and the coefficient by an element of a {\em
theory\/} which is, in an appropriate sense, dual to the operad. 
This
explains why we will construct {\it simultaneously\/} both the coloring and
the coefficient.

As the second step we offer our toy model -- the tree complex
(Section~\ref{Kolin}),
which is very easy to understand. The general case is basically the same,
only more technical, as we take into the account the symmetric group
action and the cyclic structure. 

As an application of our approach, we will try 
in Section~\ref{5} to elucidate the
construction of a cycle out of an infinitesimal deformation of a cyclic
algebra, as in~\cite{penkava-schwarz:TAMS95,penkava:infinity}.

Let us finish the introduction with the following speculations. 
The `graph complex'
for a (normal, non-cyclic) operad is the tree complex (= disguised 
bar construction), the `graph complex' for a cyclic operad is
Kontsevich's graph complex (a special case of the 
Feynman transform),
while the `graph complex' of a modular operad is the 
Feynman transform.
We believe that the construction of the classes can 
be somehow carried
over also to the last case. We will see that, for a 
Koszul (non-cyclic) 
operad, the class in the tree complex can be very 
explicitly described.
What is the homological property (if there is such a 
property) of cyclic
and/or modular operads which would made such a 
calculation possible also in
the two remaining cases?   

We assume that all algebraic objects are defined over a fixed field
$\bk$. For a vector space $V$, let $\dual V$ denote its dual, $\dual V :=
\Hom(V,\bk)$. For a natural number $k$, $\SS_k$ (resp.~$\CC_k$) denotes
the symmetric group on $k$ elements (resp.~the cyclic group of order $k$).
All calculations are made only up to signs and degrees.

\section{Warming up}
\label{Kolin}

In this section we describe the toy model of our construction -- 
the tree
complex of an operad, which is, in fact, a bar 
construction in a disguise.
The construction is easy enough not to frustrate 
the reader by unnecessary
details, but it illustrates well all the basic tricks -- 
the formulation for a
general operad (\S\ref{autumn trees}), 
the construction of an universal cycle (\S\ref{universal cycle tree}) and
the evaluation at a concrete algebra (\S\ref{the evaluation}). 
Moreover, in some lucky cases
the homology class of the universal cycle can be explicitly
described (\S\ref{Koszul case}).

\begin{odstavec}
{\rm%
{\it Autumn (colored) trees.}\label{autumn trees}
Let $\P$ be an operad and $\T$ a rooted tree. 
We say that $\T$ is $\P$-colored, if
each vertex $v$ of $\T$ with $k$ input edges 
is `colored' by an element of
$\P(k)$. Let $\ColoredTrees{\P}in$ be the 
set of all $\P$-colored trees with
$n$ input edges and $i$ inner edges.

Let $\T \in \ColoredTrees{\P}in$ and let $e$ be an inner edge of $\T$,
joining vertices $v'$ and $v''$. Define
$\partial_e(\T) \in \ColoredTrees{\P}{i-1}n$ 
to be the colored tree which is,
as a tree, obtained by the collapsing of the edge $e$, while 
the coloring of the
resulting new vertex is the obvious composition of the
corresponding colorings at $v'$ and $v''$, as indicated on 
Figure~\ref{picture1}.
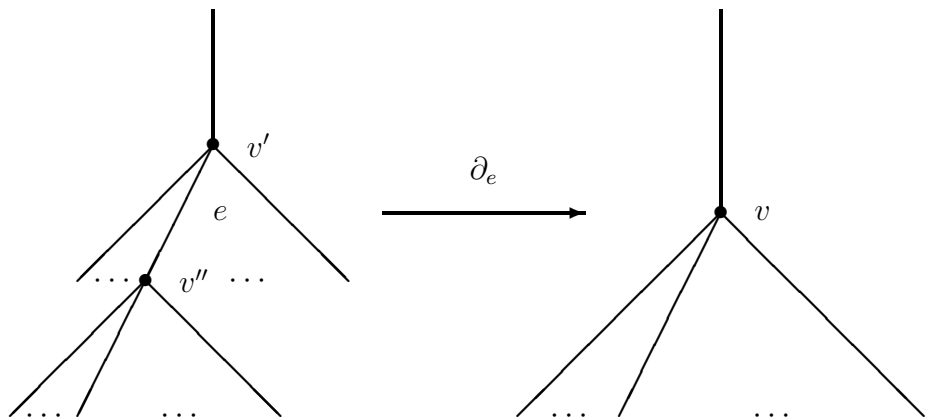
\begin{figure}
\setlength{\unitlength}{.9cm}
\begin{center}
\begin{picture}(14,6)(0,0)
\thicklines

\put(3,4){\line(0,1){2}\line(-1,-1){2}}
\put(3,4){\line(-1,-2){1}}
\put(3,4){\line(1,-1){2}}
\put(3,4){\makebox(0,0){$\bullet$}}
\put(3.5,4){\makebox(0,0)[l]{$v'$}}
\put(1.5,2){\makebox(0,0){$\cdots$}}
\put(3.5,2){\makebox(0,0){$\cdots$}}

\put(2,2){\line(-1,-1){2}}
\put(2,2){\line(-1,-2){1}}
\put(2,2){\line(1,-1){2}}
\put(2,2){\makebox(0,0){$\bullet$}}
\put(2.5,2){\makebox(0,0)[l]{$v''$}}
\put(0.5,0){\makebox(0,0){$\cdots$}}
\put(2.5,0){\makebox(0,0){$\cdots$}}

\put(3,3){\makebox(0,0)[l]{$e$}}

\put(10.5,3){\line(0,1){3}\line(-1,-1){3}}
\put(10.5,3){\line(-1,-2){1.5}}
\put(10.5,3){\line(1,-1){3}}
\put(10.5,3){\makebox(0,0){$\bullet$}}
\put(11,3){\makebox(0,0)[l]{$v$}}
\put(8.25,0){\makebox(0,0){$\cdots$}}
\put(11.25,0){\makebox(0,0){$\cdots$}}

\put(5.5,3){\vector(1,0){3}}
\put(7,3.5){\makebox(0,0)[b]{$\partial_e$}}

\end{picture}
\end{center}
\caption{
\label{picture1}
If the vertex $v'$ is colored by $c' \in {\cal P}(k)$ and 
the vertex $v''$ is colored by $c'' \in {\cal P}(l)$, 
then the resulting
vertex $v$ is colored by $c := c'(1,\ldots,c'',\ldots,1) 
\in {\cal
P}(k+l-1)$.
}
\end{figure}
The differential $\partial : \ColoredTrees{\P}in \to 
\ColoredTrees{\P}{i-1}n$
is defined by $\partial(\T) := \sum \pm \partial_e(\T)$, 
where the summation
is taken over all inner edges of the tree $\T$. 
The condition $\partial^2 =
0$ is an easy consequence of the axioms of an operad. 
The complex
\begin{equation}
\label{bar}
\ColoredTrees{\P}0n \stackrel{\partial}{\longleftarrow}
\ColoredTrees{\P}1n \stackrel{\partial}{\longleftarrow}
\cdots
\stackrel{\partial}{\longleftarrow}
\ColoredTrees{\P}{n-3}n \stackrel{\partial}{\longleftarrow}
\ColoredTrees{\P}{n-2}n
\end{equation}
is basically the bar construction $(B(\P),d_B)$ 
on the operad $\P$ (see~\cite{ginzburg-kapranov:DMJ94}), with
the opposite grading.
}\end{odstavec}

\begin{odstavec}
{\rm%
{\it The universal cycle.}\label{universal cycle tree}
Let $\cobar \P = ({\cal F}(\dual\P),\partial_\Omega)$ 
be the cobar construction on the dual cooperad $\dual\P$ and 
let $\Q$ be the operad for $\cobar\P$-algebras with  
trivial differential 
(Appendix~\ref{algebras over cobar}). 
For any $n\geq 1$ and $0 \leq i \leq n-2$ there exists a
universal cycle $\xi_i(n) \in C_i(\ColoredTrees{\P}*n; \Q(n)) :=
\ColoredTrees{\P}in \ot \Q(n)$. It is constructed as follows.

To any vertex $v$ of $\T$ with $k$ input 
edges attach the `canonical element'
$\eta \in \dual \P(k)\ot \P(k)$ corresponding, by duality, 
to the identity map $\P(k) \to \P(k)$.  
Now decorate the vertices of $\T$
with these canonical elements and interpret 
the `$\P$'-part as a coloring of
the vertex, and the `$\dual\P$'-part as an element of $\Q$ 
under the canonical
monomorphism $\dual\P \hookrightarrow \Q$ 
(Appendix~\ref{algebras over cobar}). Composing the
`$\dual\P$'-parts as indicated by the tree, 
we get an element of $\Q(n)$, where
$n$ is the number of input edges of $\T$. 
The summation over all such trees
gives $\xi_i(n)$. It follows easily that is satisfies
$\partial (\xi_i(n))=0$, thus the construction 
gives rise to an element
\[
[(\xi_i(n)] \in H_i(\ColoredTrees{\P}*n;\Q(n)).
\]
}\end{odstavec}

\begin{odstavec}
{\rm%
{\it The evaluation.}\label{the evaluation}
If $A : \Q \to \End(V)$ is a $\Q$-algebra, 
we may evaluate at $A$ to obtain
a cycle $A(\xi_i(n)) \in C_i(\ColoredTrees{\P}*n; 
\Hom(\otexp Vn,V ))$ which in
turn gives an element $A[(\xi_i(n)] \in 
H_i(\ColoredTrees{\P}*n;\Hom(\otexp
Vn,V ) )$.
}\end{odstavec}

\begin{odstavec}
{\rm%
{\it The Koszul case.}\label{Koszul case}
The reader familiar with the theory of Koszul 
operads may find interesting the
following explicit description of the element $[(\xi_i(n)]$.
If the operad $\P$ is 
Koszul~\cite{ginzburg-kapranov:DMJ94}, 
then, by the very definition of the
Koszulness and the universal coefficient formula,
\[
H_i(\ColoredTrees{\P}*n;\Q(n)) = \left\{
\begin{array}{ll}
\dual{\P^!}(n) \ot \Q(n), \mbox{ for $i=n-2$,}
\\
0, \mbox{ otherwise,}
\end{array}
\right.
\]
where $\P^!$ is the Koszul dual of the operad $\P$. 
There exists a natural
inclusion $\P^!(n) \hookrightarrow \Q(n)$ and $[(\xi_i(n)]$ 
is the
image of the canonical element of $\dual{\P^!}(n) \ot \P^!(n)$ 
under the induced
inclusion $\dual{\P^!}(n) \ot \P^!(n) \hookrightarrow \dual{\P^!}(n) 
\ot \Q(n)$.
}\end{odstavec}

\section{Graph complexes}
\label{3ii}

As we have already observed, graph complexes are 
special cases of the Feynman
transform ${\sf F}$ introduced by E.~Getzler and M.~Kapranov 
in~\cite{getzler-kapranov:preprint94}, so we could just say that
\begin{equation}
\label{Feynman}
{\cal G}^\P(n) := \bigoplus_{g\geq 0} {\sf F}\P(g,n),
\end{equation}
where ${\cal G}^\P(n)$ is the graph complex we are 
going to use, the natural
number $n$ denotes the number of external 
edges and $g$ refers to the
`genus'. We would like, however, to consider separately also the
`nonsymmetric' variant of the construction. 
Again, because a nonsymmetric
operad can be considered (after tensoring with 
the regular representation of
the symmetric group) as a symmetric operad, 
definition~(\ref{Feynman}) would
apply to this case as well, but this approach 
would obscure the {\it ribbon\/}
structure of the underlying graph. We also 
need a notation, that is why we
decided to include an explicit definition here.

\begin{odstavec}
{\rm%
{\it Symmetric vs.~\raise2pt\hbox{n}\hskip1pt\raise-1pt\hbox{o}%
\hskip-1pt\raise.5pt\hbox{n}s\raise2.5pt\hbox{y}m\raise-1pt\hbox{m}%
e\hskip2pt\raise-.5pt\hbox{t}r\hskip-2pt i\raise-3pt\hbox{c}.\/}
\label{dichotomy}
We distinguish two cases -- the nonsymmetric case and the symmetric case.
In the nonsymmetric case we work with ribbon graphs and 
nonsymmetric cyclic
operads (\S\ref{cyclic operads}), 
while in the symmetric case we work with the ordinary
cyclic operads (in the sense of~\cite{getzler-kapranov:cyclic}) 
and ordinary graphs. 

The conceptual explanation of this dichotomy is the following.
The vertices of our graphs are colored by elements of an operad. A
`color' of a vertex $v$ must behave well under the group of local
symmetries of the graph at $v$. How does this group look? 
For a general
graph, it is the group $\SS_{k+1}$ 
permuting the (half)edges at $v$, $k+1$ being the 
number of these edges. This
means that the `color' at $v$ must admit a 
$\SS_{k+1}$-symmetry, and we
necessary 
arrive at the notion of a cyclic operad. 
In a ribbon graph, the set of
(half)edges at $v$ has a preferred cyclic order. 
The group $\SS_{k+1}$ is the
semidirect product of $\SS_k$ and the cyclic group 
$\CC_{k+1}$, and the cyclic
order of the edges fixes the $\SS_k$-part, thus t
he `color' 
at $v$ must admit a
$\CC_{k+1}$-symmetry. The corresponding notion is 
that of a nonsymmetric
cyclic operad, see~\S\S\ref{cyclic operads}, 
\ref{colored graphs} and
\ref{contracting an edge} for details. 
}
\end{odstavec}

\begin{odstavec}
{\rm%
{\it Graphs.\/}\label{graphs}
As usual, a graph consists of edges and vertices. We suppose
that all vertices are at least trivalent. 
Let $\vert{\G}$ denote the set of vertices of the graph
$\G$ and let  $\edg{\G}$ be the set of edges 
of $\G$. 
For a vertex $v\in \vert {\G}$ let $\edg v$ denote the
set of half edges (there may be loops in the graph!) 
at the vertex $v$. A ribbon graph is a graph such that 
a cyclic order on the set $\edg v$ is given, 
for any $v\in \vert {\G}$.
We allow our graph $\G$ to have also some external
edges; we denote the set of all these external edges by $\leg {\G}$.
Let $\Graphs in$ denote the set of graphs with $n$ 
external edges and $i$
internal edges attached to two distinct vertices.
}\end{odstavec}

\begin{odstavec}
{\rm%
{\it Cyclic operads.\/}\label{cyclic operads} 
Following~\cite{getzler-kapranov:cyclic}, 
a cyclic operad is an operad $\P$ such that the usual
action of the symmetric group $\SS_n$ on $\P(n)$ 
is extended to an action of
the symmetric group $\SS_{n+1}$. This extension 
has, of course, to satisfy
appropriate axioms. There exists 
an (almost) obvious nonsymmetric
version where each $\P(n)$ has an action of the 
cyclic group $\CC_{n+1}$. We call these objects nonsymmetric
cyclic operads.
}\end{odstavec}

\begin{odstavec}
{\rm%
{\it $\P$-colored graphs.\/}\label{colored graphs} 
If $I$ is a cyclically ordered
set of $n+1$ elements, then the cyclic group $\CC_{n+1}$ acts on
the set of all cyclic-%
order preserving maps $f: \{0,1,\ldots,n\} \to I$. If $V$ is a
$\CC_{n+1}$ space we put, as in~\cite{ginzburg-kapranov:DMJ94},
\[
V((I)):= \left(
\bigoplus_{f: \{0,1,\ldots,n\} \to I} V
\right)_{\CC_{n+1}}
\ \mbox{(the set of coinvariants).}
\]
By a $\P$-colored graph we mean a graph 
$\G$ such that each vertex $v$ is `colored' by an element of
$\P((\edg v)$, where the cyclic order on $\edg v$ i
s given by the ribbon
structure of the graph. The symmetric variant of 
this definition is obvious.
Denote by $\ColoredGraphs{\P}in$ the set of 
$\P$-colored graphs with $n$
external and $i$ internal edges attached to two distinct vertices.
}\end{odstavec}

\begin{odstavec}
{\rm%
{\it Contracting an edge.\/}\label{contracting an edge}
Let us discuss the nonsymmetric case first.
Let $\G \in \ColoredGraphs{\P}in$ be a $\P$-colored graph and  let
$e \in \edg {\G}$ be an edge attached
to two distinct vertices. 
Define the $\P$-colored graph $\partial_e(\G)\in 
\ColoredGraphs{\P}{i-1}n$
as follows. As a graph it
coincides with the graph $\G/e$ obtained by collapsing 
out the edge $e$ from
$\G$, with the induced cyclic order on the 
resulting vertex $v$. 

Before going further, we need some notation. 
Let $e$ join (distinct)
vertices $v'$ and $v''$ and let $\edg{v'} = 
\set{e,e'_1,e'_2,\ldots,e'_k}$, 
$\edg {v''} = \set{e,e''_1,e''_2,\ldots,e''_l}$ (in
this cyclic order). This means that $\edg v =
\set{e''_l,e'_1,\ldots,e'_k,e''_1,\ldots,e''_{l-1}}$ 
(in this cyclic 
order), see Figure~\ref{picture2}.
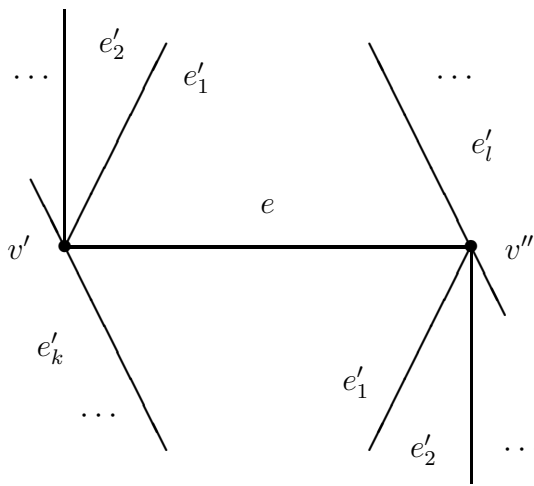
\begin{figure}
\setlength{\unitlength}{.9cm}
\begin{center}
\begin{picture}(8,6)(0,1)
\thicklines

\put(1,4){\line(1,2){1.5}}
\put(1,4){\line(0,1){3.5}}
\put(1,4){\line(1,-2){1.5}}
\put(.5,5){\line(1,-2){.5}}
\put(1,4){\makebox(0,0){$\bullet$}}

\put(.5,4){\makebox(0,0)[r]{$v'$}}
\put(.5,6.5){\makebox(0,0){$\cdots$}}
\put(1.5,1.5){\makebox(0,0){$\cdots$}}

\put(2.75,6.5){\makebox(0,0)[l]{$e'_1$}}
\put(1.5,7){\makebox(0,0)[l]{$e'_2$}}
\put(1,2.5){\makebox(0,0)[r]{$e'_k$}}

\put(1,4){\line(1,0){6}}
\put(4,4.5){\makebox(0,0)[b]{$e$}}

\put(7,4){\line(-1,-2){1.5}}
\put(7,4){\line(0,-1){3.5}}
\put(7,4){\line(-1,2){1.5}}
\put(7,4){\line(1,-2){.5}}
\put(7,4){\makebox(0,0){$\bullet$}}

\put(7.5,4){\makebox(0,0)[l]{$v''$}}
\put(7.75,1){\makebox(0,0){$\cdots$}}
\put(6.75,6.5){\makebox(0,0){$\cdots$}}

\put(7,5.5){\makebox(0,0)[l]{$e'_l$}}
\put(5.5,2){\makebox(0,0)[r]{$e'_1$}}
\put(6.5,1){\makebox(0,0)[r]{$e'_2$}}

\end{picture}
\end{center}
\caption{An edge $e$ joining the vertices $v'$ and $v''$.
\label{picture2}
}
\end{figure}
Let $c' \in \P((\edg{v'}))$ and $c'' \in \P((\edg{v''}))$
be the colorings of the vertices $v'$ and $v''$, respectively.

Define $f' : \set {0,1,\ldots,k} \to \edg {v'}$ by  $f'(0) =e$, and
$f'(i) = e'_i$ for $1\leq l \leq  k$, and let $r' 
\in \P(k+1)_{f'}$ be a
representative for $c'$. Similarly, let 
$f'' : \set {0,1,\ldots,l} \to \edg {v''}$ be given 
by $f''(0)=e''_l$, $f''(1)
= e$ and $f''(j) = e''_{j-1}$ for $2\leq j \leq l$, and let 
$r'' \in \P(l+1)_{f''}$ be a
representative for $c''$. We then define the coloring $c$
of $v$ to be the equivalence class of $r''(r',1,\ldots,1)$  
in $\P(k+l-1)_f$, where
$f(0) = e''_l$,
$f(i) = e'_i$, for $1\leq i\leq k$, and $f(k+j) = e''_j$, 
for $1\leq j\leq
l-1$. It follows from the cyclicity of the operad $\P$ 
that the coloring $c$
does not depend on the particular choice of the 
representatives $c'$ and
$c''$. The symmetric
case is even easier, because 
we do not need not to pay any attention to the cyclic
order of the edges. 
}\end{odstavec}

\begin{odstavec}
{\rm%
{\it Graph complex.\/}\label{graph complex}
The differential 
$\partial :\ColoredGraphs{\P}in \to 
\ColoredGraphs{\P}{i-1}n$ is defined by
$\partial(\G) = \sum \pm \partial_e(\G)$, 
the sum being taken over all edges
attached to two distant vertices. The condition
$\partial^2 = 0$ follows from the fact that 
$\P$ is an operad. We call 
$\ColoredGraphs{\P}*n = (\ColoredGraphs{\P}*n,\partial)$ 
the $\P$-graph
complex.
For $n=0$ we write
simply ${\cal G}^{\P}_*$ instead of $\ColoredGraphs{\P}*0$.
}\end{odstavec}

\begin{odstavec}
{\rm%
{\it Examples.\/}\label{examples}
The nonsymmetric 
operad $\Ass$ for associative algebras is a 
nonsymmetric cyclic
operad. Because $\Ass(n) = \bk$ for each $n\geq 1$,
the coloring contains no information and 
the complex ${\cal G}^{\Ass}_* = 
({\cal G}^{\Ass}_*,\partial)$ is the
complex of ribbon graphs introduced in~\cite{kontsevich:91}. 

Similarly, the symmetric cyclic operad $\Comm$ for
commutative algebras has $\Comm(n) = \bk$ and 
the complex 
${\cal G}^{\Comm}_* = ({\cal G}^{\Comm}_*,\partial)$ is the
complex of (all) graphs considered 
in~\cite{kontsevich:91,penkava:infinity}.

For the symmetric cyclic operad $\Lie$, the graph complex 
${\cal G}^{\Lie}_* = ({\cal G}^{\Lie}_*,\partial)$ 
was constructed 
by~M.~Kontsevich in~\cite{kontsevich:93}. 
We may imagine a $\Lie$-graph
as a graph whose vertices are colored by $(k-1)!$ `colors' 
representing a
basis of $\Lie(k)$, where $k+1$ is the number of 
edges at the vertex.
}\end{odstavec}

\section{The cycle}
\label{4}

\begin{odstavec}
{\rm%
{\it The `state sum'.\/}\label{state sum}
Let $\theor$ be the theory describing cyclic 
$\cobar{\P}$-algebras with
trivial differential (Appendix~\ref{algebras over cobar}). 
The universal cycle
$\xi_i(n) \in C_i(\ColoredGraphs{\P}*n;\theor(n,0))$ 
is defined as follows.
Decorate each vertex $v$ of $\G \in \Graphs in$ 
with the `canonical element'
of $\dual\P(k) \ot \P(k)$, $k = 
\mbox{ord}\{\edg v\}$. 
As in the case of trees,
we interpret the `$\P$'-part as the 
coloring of the vertex 
$v$ and the `$\dual\P$'-part as
an element of $\theor(k+1,0)$, under 
the canonical monomorphism $\dual\P(*)
\hookrightarrow \theor(*+1,0)$ 
(Appendix~\ref{cyclic algebras}). 
Now compose these `$\dual\P$'-parts as
elements of $\theor$, using $\nu \in \theor(0,2)$ 
as a `propagator' along the
edges of $\G$. This gives an element of $\theor(n,0)$. 
The requisite
$\xi_i(n)$ is then the summation over all 
graphs $\G \in \Graphs in$.
}\end{odstavec}

The central statement is the following proposition.
\begin{proposition}
\label{central}
The chain $\xi_i(n) \in C_i(\ColoredGraphs{\P}*n; 
\theor(n,0))$ is a cycle,
$\partial(\xi_i(n)) =0$. This means that it 
determines a homology class
\[
[\xi_i(n)] \in H_i(\ColoredGraphs{\P}*n; \theor(n,0)).
\]
\end{proposition}

\begin{odstavec}
{\rm%
{\it The cycle defined by a cyclic algebra.\/}\label{cycle}
Let $B = (V,A,h,\nu)$ be a cyclic $\cobar\P$-algebra as 
in~\S\ref{cyclic algebras}, 
i.e.~a map $B: \theor \to
\END(V)$ of theories. 
The evaluation at $B$ gives
a cycle $B(\xi_i(n)) \in 
C_i(\ColoredGraphs{\P}*n;\Hom(\otexp Vn,\bk))$ which
in turn defines the class
\[
B([\xi_i(n)]) \in H_i(\ColoredGraphs{\P}*n;\Hom(\otexp Vn,\bk)).
\]
Extremely important is the case $n=0$ when $\Hom(\otexp
Vn,\bk) = \bk$. We get elements
\[
B(\xi_i) \in {\cal G}^{\P}_i \mbox{ and }
B([\xi_i]) \in H_i({\cal G}^{\P}_*).
\]
}\end{odstavec}

\begin{odstavec}
{\rm%
{\it Examples.\/}\label{examples2}
If $\P = \Ass$, the nonsymmetric operad 
for associative algebras,
Proposition~\ref{central} gives the `universal element' 
\[
[\xi] \in H({\cal M}_{g,n}; \theor(0,0)).
\]
A cyclic $\cobar{\Ass}$-algebra is an 
$A_\infty$-algebra with nondegenerate
invariant scalar product. The evaluation at 
such an algebra gives the classes
in $H({\cal M}_{g,n}; \bk)$ constructed by M.~Kontsevich 
in~\cite{kontsevich:91}.

If $\P = \Comm$, the symmetric operad for 
associative commutative algebras,
then a cyclic  $\cobar{\Comm}$-algebra 
is an $L_\infty$ (or strong homotopy Lie) algebra 
with nondegenerate
invariant scalar product. Evaluation at $[\xi_i]$ 
then gives the classes
constructed in~\cite{penkava:infinity}.
}\end{odstavec}

\section{Lacunar graphs}
\label{5}

The construction of M.~Kontsevich of 
an element in $H^*({\cal M}_{g,n};
{\bf C})$ requires a cyclic $A_\infty$-algebra; 
we already know that this
element is an evaluation at the universal element of
Proposition~\ref{central}. In~\cite{penkava-schwarz:TAMS95}, 
M.~Penkava and 
A.~Schwarz showed that, if we
restrict to a suitable subcomplex of the graph complex, 
we may construct a
similar invariant out of an infinitesimal deformation of a cyclic
associative algebra. M.~Penkava then 
generalized in~\cite{penkava:infinity} 
this construction to the $L_\infty$ (strong homotopy~Lie)
algebra case. Here we explain an almost obvious generalization of 
these results to algebras over an arbitrary cyclic operad
and give also a conceptual explanation of the construction.

Since we work with deformations, we need an 
independent variable $t$. We work
over the extended coefficient ring $\bk [t]$, the ring of
polynomials in $t$. We will
use the notation $\P[t]$ for $\P \ot_{\bk}\bk [t]$,~etc. 

Fix $k\geq 3$ and consider (for a fixed $i$) the subspace
$\ColoredGraphsPrime{\P}in$ of $\ColoredGraphs{\P[t]}in =
\ColoredGraphs{\P}in[t]$ consisting of 
$\P[t]$-colored graphs such
that
\begin{itemize}
\item[(i)]
all vertices are trivalent except exactly one 
which is $(k+1)$-valent,
\item[(ii)]
trivalent vertices are colored by elements of 
$\P(2) =\P(2)\cdot t^0 
\subset\P(2) [t]$ and
\item[(iii)]
the $(k+1)$-valent vertex is colored by an 
element of $\P(k)\cdot t^1 
\subset\P(k) [t]$.  
\end{itemize} 
Now construct the universal cycle $\xi_i(n) 
\in C_i(\ColoredGraphs{\P[t]}*n,
\theor(n,0)[t])$ and restrict it to a cycle 
$\xi'_i(n) \in C_i(\ColoredGraphsPrime{\P}*n,
\theor(n,0)[t])$. The crucial observation is 
that $\partial(\xi'_i(n))$
consists of at
most linear terms in $t$. So, to evaluate at $\xi'_i(n)$, 
a map $A: \theor[t] \to \END(V)[t]$ which is 
a map of theories 
{\it \underline{modulo}\/} $\underline{t}^2$
is enough!

A moment's reflection shows that an infinitesimal 
deformation of a cyclic 
$\P^!$-algebra into a cyclic $\cobar{\P}$-algebra 
gives such a map. We may
conclude this paragraph by observing that 
infinitesimal deformations of cyclic
algebras over a Koszul cyclic operad
are governed by the cyclic cohomology which was 
constructed, for a general
cyclic operad, in~\cite{ginzburg-kapranov:DMJ94}.  

\section{Appendix}
\label{6}

\begin{odstavec}
{\rm%
{\it Algebras over $\cobar\P$.\/}\label{algebras over cobar}
Let $\P$ be an  operad and let $\cobar{\P}$ be the cobar
construction~\cite{ginzburg-kapranov:DMJ94} on the dual cooperad 
$\dual\P$, $\dual \P = \{\dual \P(n)\}_{n\geq 1}$, 
with $\dual \P(n) =
\dual{\P(n)}= \Hom(\P(n),\bk)$. 
This means that $\cobar{\P} =
(\F(\dual\P), \partial_\Omega)$, where  $\F(\dual\P)$ is
the free operad on the collection $\dual\P$ and the differential
$\partial_\Omega$ is induced by the structure maps of 
the cooperad $\dual\P$.

As usual, an algebra over $\cobar\P$ is a differential 
vector space $(V,d_V)$
together with a map $A: \cobar\P \to \End(V,d_V)$ of 
differential operads;
here $\End(V,d_V)$ is the endomorphism operad of $(V,d_V)$. 
As we work with
graphs having at least trivalent vertices, we consider 
only the case
$d_V = 0$. Such algebras can be described as algebras 
over the operad $\Q =
\F(\dual\P)/(\partial_\Omega(\dual\P))$ (= the free 
operad on $\dual\P$ modulo the
ideal generated by $\partial_\Omega(\dual\P)$). As an 
easy consequence of the
quadraticity of the differential $\partial_\Omega$ we 
see that the canonical
projection $\dual\P \to \F(\dual\P)$ induces a 
monomorphism $\dual \P \hookrightarrow
\Q$ of collections.

}\end{odstavec}

\begin{odstavec}
{\rm%
{\it Cyclic algebras.\/}\label{cyclic algebras}
Let $A$ be an $\cobar\P$-algebra as above. We say that $A$
is cyclic if there exists a symmetric bilinear product 
$h=\scalar --$ on $V$
such that
\begin{equation}
\label{scalar}
\scalar {A(\phi)(x_1,\ldots,x_n)}{x_{n+1}} = 
\scalar {x_{1}}{A(\phi)(x_2,\ldots,x_{n+1})},
\end{equation}
for any $\phi \in \dual\P(n)$, and $x_1,\ldots,x_{n+1}\in V$.

Thus, a cyclic algebra is an object of the form 
$B=(V,A,h,\nu)$, where $A: \Q
\to \End(V)$ is a $\Q$-algebra structure on $V$, 
$h: V\ot V \to \bk$ is a
scalar product and $\nu : \bk \to V\ot V$ is 
the `inverse matrix' of $h$ 
in the sense that
\[
(h\ot \id)\circ (\id \ot \nu) = (\id \ot h)\circ 
(\nu \ot \id) = \id.
\]
Such objects form an equationally given category which 
is not algebraic,
i.e.~it can not be described as a category of 
algebras over an operad, but
rather as a category of algebras over a 
$\bk$-linear {\small PROP}. 
In~\cite{markl:JPAA96,markl:ws93} 
we used the name `theory' for a $\bk$-linear {\small PROP}.
Although this terminology is obviously
not the best one, we will use this name here. 
We refer to~\cite{markl:JPAA96,markl:ws93} 
for a very thorough introduction to
$\bk$-linear {\small PROP}s. 

Just recall that a theory is a sequence $\theor
= \{ \theor(m,n);\ m,n\geq 0\}$ of $\bk$-linear spaces, 
each $\theor(m,n)$ encoding operations with
$m$ inputs and $n$ outputs, as the 
$\P(n)$-part of an operad $\P$ 
encodes operations with $n$
inputs and just one output. 
For any vector space $V$ there exists the
`endomorphism theory' $\END(V)$, with $\END(V)(m,n) = 
\Hom(\otexp Vm, \otexp
Vn)$. An algebra over $\theor$ is then a map $B: 
\theor \to \END(V)$ of
theories. 

Let $\theor$ be the theory describing cyclic $\Q$-algebras. 
The theory
$\theor$ is, in a well defined sense, 
generated by the operad $\Q$ and by the
elements $h\in \theor(2,0)$ and $\nu \in \theor(0,2)$.
The correspondence $q\mapsto h(q,\id)$ 
defines a map $\Q(k) \to
\theor(k+1,0)$. The composition of this 
map with the inclusion $\dual\P(k) \to
\Q(k)$ gives the inclusion $\dual \P(k) 
\hookrightarrow \theor(k+1,0)$. The
elements in the
image of this map are, due to~(\ref{scalar}), 
invariant under the action of
the cyclic group.  
}\end{odstavec}


\catcode`\@=11
\noindent

\vskip3mm
\catcode`\@=11
\noindent
M.~M.: Mathematical Institute of the Academy, 
\v Zitn\'a 25, 115 67
Praha 1, Czech Republic,\hfill\break\noindent
\hphantom{M.~M.:} email: 
{\tt markl@math.cas.cz}\hfill\break\noindent

\end{document}